# NEW APPROXIMATIONS FOR THE AREA OF THE MANDELBROT SET

DANIEL BITTNER, LONG CHEONG, DANTE GATES, AND HIEU D. NGUYEN

ABSTRACT. Due to its fractal nature, much about the area of the Mandelbrot set $M$ remains to be understood. While a series formula has been derived by Ewing and Schober to calculate the area of $M$ by considering its complement inside the Riemann sphere, to date the exact value of this area remains unknown. This paper presents new improved upper bounds for the area based on a parallel computing algorithm and for the 2-adic valuation of the series coefficients in terms of the sum-of-digits function.

## 1. INTRODUCTION

The Mandelbrot set (hereafter $M$) is defined as the set of complex numbers $c \in \mathbb{C}$ such that the sequence $\{z_n\}$ defined by the recursion
$$z_n = z_{n-1}^2 + c \tag{1}$$
with initial value $z_0 = 0$ remains bounded for all $n \geq 0$. Douady and Hubbard [6] proved that $M$ is connected and Shishikura [18] proved that $M$ has fractal boundary of Hausdorff dimension 2. However, it is unknown whether the boundary has positive Lebesgue measure.

Ewing and Schober [9] derived a series formula for the area of $M$ by considering its complement, $\tilde{M}$, inside the Riemann sphere $\overline{\mathbb{C}} = \mathbb{C} \cup \{\infty\}$, i.e. $\tilde{M} = \overline{\mathbb{C}} - M$. It is known that $\tilde{M}$ is simply connected with mapping radius 1 ([6]). In other words, there exists an analytic homeomorphism
$$\psi(z) = z + \sum_{m=0}^{\infty} b_m z^{-m} \tag{2}$$
which maps the domain $\Delta = \{z : 1 < |z| \leq \infty\} \subset \overline{\mathbb{C}}$ onto $\tilde{M}$. It follows from the classic result of Gronwall [11] that the area of the Mandelbrot set $M = \overline{\mathbb{C}} - \tilde{M}$ is given by
$$A = \pi \left[ 1 - \sum_{m=1}^{\infty} m|b_m|^2 \right] \tag{3}$$

The arithmetic properties of the coefficients $b_m$ have been studied in depth, first by Jungreis [12], then independently by Levin [14, 15], Bielefeld, Fisher, and Haeseler [3], Ewing and Schober [8, 9], and more recently by Shimauchi [17]. We note that the results of Levin [14, 15] and Shimauchi [17] hold for Multibrot sets defined by generalizing (1) to higher-order recurrences.

There are three approaches to calculating the coefficients $b_m$. The first approach involves expressing $b_m$ as a contour integral, found independently by Levin [14] and by Ewing and Schober [8]:
$$b_m = -\frac{1}{2\pi m i} \int_{|z|=R} p_n(z)^{m/2^n} dz \tag{4}$$
where $1 \leq m \leq 2^{n+1} - 3$ and $R$ is chosen sufficiently large. The polynomials $p_n(w)$ in (4) are defined recursively by
$$\begin{aligned} p_0(w) &= w \\ p_n(w) &= p_{n-1}^2(w) + w. \end{aligned} \tag{5}$$
Ewing and Schober proved [8] that the polynomials $p_n(w)$ are Faber polynomials of degree $2^n$ for $M$, i.e. $p_n(\psi(z)) = z^{2^n} + o(1)$ as $z \to \infty$, a fact that they used to prove (4). Jungreis [12] proved earlier that

---

$b_{2^n+1} = 0$ for $n \geq 1$ (see also [3, 8, 14]). Bielefeld, Fisher, and Haeseler [3] proved that no constants $\epsilon$ and $K$ exist so that $|b_m| < K/m^{1+\epsilon}$ for all $m$.

The second approach to calculating $b_m$, due to Bielefeld, Fisher, and Haeseler [3], involves substituting (2) into (5) to obtain
$$p_n(\psi(z)) = p_{n-1}^2(\psi(z)) + \psi(z) = z^{2^n} + o(1)$$
and then equating coefficients to recursively solve for $b_m$. In this paper, we follow a variation of this approach, due to Ewing and Schober [9], by expanding $p_n(\psi(z))$ in the form
$$p_n(\psi(z)) = \sum_{m=0}^{\infty} \beta_{n,m} z^{2^n - m} \tag{6}$$
where $b_m = \beta_{0,m+1}$. It follows that $\beta_{n,m} = 0$ for $n \geq 1$ and $1 \leq m \leq 2^n$. Moreover, this range of zero values can be extended to $1 \leq m \leq 2^{n+1} - 2$ because of the recursion
$$\beta_{n,m} = 2\beta_{n-1,m} + \sum_{k=1}^{m-1} \beta_{n-1,k} \beta_{m-1,m-k}, \tag{7}$$
which can be derived by substituting (6) into (5) and equating coefficients. Formula (7) can then be manipulated to obtain the following backward recursion formula [9]:
$$\beta_{n,m} = \frac{1}{2}\left[\beta_{n+1,m} - \sum_{k=2^{n+1}-1}^{m-2^{n+1}+1} \beta_{n,k}\beta_{n,m-k} - \beta_{0,m-2^{n+1}+1}\right] \tag{8}$$
where $\beta_{n,0} = 1$ and $\beta_{0,m} = b_{m-1}$ for $m \geq 1$.

No explicit formula is known for $\beta_{n,m}$ (nor for $b_m$), except those at certain positions. However, it is clear from (8) that $\beta_{n,m}$ is rational and that its denominator equals a power of 2 when expressed in lowest terms. In [9] Ewing and Schober established the following upper bound on its 2-adic valuation.

**Theorem 1** (Ewing-Schober [9]). *Let $n \in \mathbb{N}$ and $m$ be a positive integer. Then $2^{2m+3-2^{n+2}} \beta_{n,m}$ is an integer, i.e.*
$$-\nu(\beta_{n,m}) \leq 2m + 3 - 2^{n+2} \tag{9}$$
*for non-zero $\beta_{n,m}$.*

Here, the 2-adic valuation $\nu(x)$ of a positive integer $x$ is defined to be the greatest integer for which $2^{\nu(x)}$ divides $x$ and if $x/y$ is a fraction in lowest terms, then we define $\nu(x/y) = \nu(x) - \nu(y)$. If $x = 0$, then we set $\nu(x) = \infty$. Observe that in the special case $b_m = \beta_{0,m+1}$, (9) reduces to
$$-\nu(b_m) \leq 2m - 1. \tag{10}$$

Zagier [3] observed earlier that
$$-\nu(b_m) \leq \nu((2m+2)!)$$
for $0 \leq m \leq 1000$. Moreover, he observed that equality holds when $m$ is odd (or zero). These results were later proven by Levin [14] and Shimauchi [17].

**Theorem 2** (Levin [14]). *If $m$ is a positive odd integer, then*
$$-\nu(b_m) = \nu((2m+2)!). \tag{11}$$

**Theorem 3** (Shimauchi [17]). *Let $m$ be a non-negative integer. Then*
$$-\nu(b_m) \leq \nu((2m+2)!). \tag{12}$$
*Moreover, equality holds precisely when $m$ is odd.*

Ewing and Schober [9] used (8) to compute the first 240,000 coefficients for $b_n$ by computer. Since
$$A \leq A_N \equiv \pi \left[1 - \sum_{m=1}^{N} m|b_m|^2\right], \tag{13}$$
their calculation of $A_{240,000} \approx 1.7274$ yielded an upper bound for the area of $M$. They were able to slightly improve their result to 1.72 by extending their computations to the first 500,000 coefficients as reported by Ewing [7]. They also calculated a crude lower bound of $7\pi/16 \approx 1.3744$ by estimating the size of the main cardioid ($3\pi/8$) and the main bulb ($\pi/16$). However, they reported a discrepancy with their approximation of 1.52 obtained by pixel counting. More recent calculations by Förstemann [10] provide an estimate of 1.50659 based on a resolution of almost 88 trillion pixels. In addition, Andreadis and Karakasidis [2] obtained an estimate of 1.5052 based on the boundary scanning method. Thus, as noted by Ewing and Schober, either the series (3) converges so slowly that the approximation $A_{500,000} \approx 1.72$ is poor or else the pixel counting method fails to account for the boundary of $M$. Recently, Buff and Chéritat [4] found Julia sets with positive area. Therefore, coupled with Shishikura's result that the boundary of $M$ has Hausdorff dimension 2, it is not far-fetched to suspect that the boundary of $M$ may have positive area.

In this paper, we report on progress in obtaining new upper bounds for $A$ and new results involving the two-dimensional sequence $\beta_{n,m}$. In particular, we were able to compute the first five million coefficients for $b_n$ by developing a parallel processing implementation of (8). This extends the calculation of the first one million coefficients by Chen, Kawahira, Li, and Yuan [5] by five-fold where they reported the upper bound $A_{1,000,000} = 1.703927$. As a result of our calculations, we obtained the new upper bound

$$A_{5,000,000} \approx 1.68288 \tag{14}$$

Moreover, we were able to improve on (9) by establishing the tighter bound (Theorem 9)

$$-\nu(\beta_{n,m}) \leq 2m - 2^{n+2} + 4 - s(n,m). \tag{15}$$

for non-zero $\beta_{n,m}$ where $s(n,m)$ is the base-2 sum-of-digits function of degree $n$ (Definition 4). In the special case $b_m = \beta_{0,m+1}$, we obtain as a corollary

$$-\nu(b_m) \leq 2(m+1) - s(0, m+1). \tag{16}$$

This is equivalent to Shimauchi's result (12) because of the relation $\nu(k!) = k - s(0,k)$ for any positive integer $k$. Observe that equality in (16) holds for all odd values of $m$, which follows from Shimauchi's result (Theorem 3), whereas (10) holds only when $m+1$ equals a power of 2.

Our new upper bound (15) is significant on two levels. First, from a computational perspective, it allows the values of $\beta_{n,m}$ to be calculated by integer arithmetic (as discussed by Ewing and Schober [9]) using less memory than (9). Such an approach would increase the accuracy in which upper bounds for the area of the $M$ are calculated over floating-point arithmetic where the values of $\beta_{n,m}$ are stored as truncated decimals. Secondly, (16) confirms Levin's work that the sum-of-digits function is a crucial ingredient in determining the exact area of $M$ by using the series formula (3).

2. Two-adic Valuation of $\beta_{n,m}$

In this section we consider the 2-adic valuation of $\beta_{n,m}$ and prove the bound (15), which is a refinement of (9). We begin by defining the sum-of-digits function and present a series of lemmas on properties of this function that will be utilized in the proof. Throughout this paper, $\mathbb{N}$ denotes the set of non-negative integers.

**Definition 4.** Let $m \in \mathbb{N}$ with base-2 expansion $m = d_L 2^L + d_{L-1} 2^{L-1} + ... + d_0 2^0$ where $d_L = 1$ and $d_i \in \{0, 1\}$ for $i < L$. We define the *base 2 sum-of-digits function $s(n,m)$ of degree $n$* by

$$s(n,m) = \sum_{i=n}^{L} d_i$$

**Lemma 5.** *Let $m, n \in \mathbb{N}$. Then $s(n,m)$ is sub-additive, i.e.*

$$s(n, l+m) \leq s(n,l) + s(n,m)$$

*for all $l \in \mathbb{N}$.*

*Proof.* We follow the proof in [16]. Let $l = c_K 2^K + c_{K-1} 2^{K-1} + ... + c_0 2^0$ and $m = d_L 2^L + d_{L-1} 2^{L-1} + ... + d_0 2^0$. Since $s(n, m + 2^i) \leq s(n, m)$ for $i < n$ and $s(n, m + 2^i) \leq s(n, m) + 1$ for $i \geq n$, it follows that

$$s(n, l + m) = s\left(n, m + \sum_{i=0}^{K} c_i 2^i\right)$$
$$\leq s\left(n, m + \sum_{i=n}^{K} c_i 2^i\right)$$
$$\leq s(n, m) + \sum_{i=n}^{K} c_i$$
$$\leq s(n, m) + s(n, l)$$

as desired. □

**Lemma 6.** *For all $m, n \in \mathbb{N}$, we have*

*a)* $0 \leq s(n, m) - s(n+1, m) \leq 1$.

*b)* $s\left(n, 2^{n+1} - 1\right) = 1$.

*c)* $s(0, m) \leq 2s\left(0, \frac{m}{2}\right) - 1$ *for positive even integers $m$.*

*Proof.* (a) We express $m$ as in Definition 4. It follows that

$$s(n, m) - s(n+1, m) = \sum_{i=n}^{L} d_i - \sum_{i=n+1}^{L} d_i$$
$$= d_n + \sum_{i=n+1}^{L} d_i - \sum_{i=n+1}^{L} d_i$$
$$= d_n$$

where $d_n$ must equal either 0 or 1. This completes the proof for part a).

b) The result follows immediately from the fact that $2^{n+1} - 1 = 2^0 + ... + 2^{n-1} + 2^n$.

c) Assume $m$ is even. Then $m$ can be expressed as

$$m = \sum_{i=r}^{L} d_i 2^i$$

for some integers $r, L$, where $r \geq 1$ by assumption. It follows that

$$\frac{m}{2} = \sum_{i=r-1}^{L-1} d_i 2^i.$$

Therefore,

$$s(0, m) = s\left(0, \frac{m}{2}\right)$$
$$= 2s\left(0, \frac{m}{2}\right) - s\left(0, \frac{m}{2}\right)$$
$$\leq 2s\left(0, \frac{m}{2}\right) - 1$$

since $s\left(0, \frac{m}{2}\right) \geq 1$. □

Next, we present a lemma regarding the convolution described in equation (8).

**Lemma 7.** *Let $m \in \mathbb{N}$ with $m \geq 2^{n+2} - 2$.*

*a) For $m$ even, we have*

$$\sum_{k=2^{n+1}-1}^{m-2^{n+1}+1} \beta_{n,k}\beta_{n,m-k} = 2\left[\sum_{k=2^{n+1}-1}^{m/2-1} \beta_{n,k}\beta_{n,m-k}\right] + \left(\beta_{n,m/2}\right)^2 \qquad (17)$$

*b) For $m$ is odd, we have*

$$\sum_{k=2^{n+1}-1}^{m-2^{n+1}+1} \beta_{n,k}\beta_{n,m-k} = 2\left[\sum_{k=2^{n+1}-1}^{(m-1)/2} \beta_{n,k}\beta_{n,m-k}\right] \qquad (18)$$

*Proof.* When $m$ is even, we have

$$\sum_{k=2^{n+1}-1}^{m-2^{n+1}+1} \beta_{n,k}\beta_{n,m-k} = \sum_{k=2^{n+1}-1}^{m/2-1} \beta_{n,k}\beta_{n,m-k} + \sum_{m/2}^{m/2} \beta_{n,k}\beta_{n,m-k} + \sum_{m/2+1}^{m-2^{n+1}+1} \beta_{n,k}\beta_{n,m-k}$$

Letting $h = m - k$, we obtain

$$\sum_{k=2^{n+1}-1}^{m-2^{n+1}+1} \beta_{n,k}\beta_{n,m-k} = \sum_{k=2^{n+1}-1}^{m/2-1} \beta_{n,k}\beta_{n,m-k} + \left(\beta_{n,m/2}\right)\left(\beta_{n,m/2}\right) + \sum_{h=m/2-1}^{2^{n+1}-1} \beta_{n,m-h}\beta_{n,h}$$

$$= 2\left[\sum_{k=2^{n+1}-1}^{m/2-1} \beta_{n,k}\beta_{n,m-k}\right] + \left(\beta_{n,m/2}\right)^2$$

This proves part a).

On the other hand, when $m$ is odd,

$$\sum_{k=2^{n+1}-1}^{m-2^{n+1}+1} \beta_{n,k}\beta_{n,m-k} = \sum_{k=2^{n+1}-1}^{(m-1)/2} \beta_{n,k}\beta_{n,m-k} + \sum_{k=(m+2)/2}^{m-2^{n+1}+1} \beta_{n,k}\beta_{n,m-k}$$

Letting $l = m - k$, we have

$$\sum_{k=2^{n+1}-1}^{m-2^{n+1}+1} \beta_{n,k}\beta_{n,m-k} = \sum_{k=2^{n+1}-1}^{(m-1)/2} \beta_{n,k}\beta_{n,m-k} + \sum_{l=(m-1)/2}^{2^{n+1}-1} \beta_{n,m-l}\beta_{n,l}$$

$$= 2\left[\sum_{k=2^{n+1}-1}^{(m-1)/2} \beta_{n,k}\beta_{n,m-k}\right]$$

This justifies part b). □

We now present one final lemma involving the right hand side of (15).

**Lemma 8.** *Let $m, n \in \mathbb{N}$ and define*

$$p(n, m) = 2m - 2^{n+2} + 4 - s(n, m). \qquad (19)$$

*Then the following inequalities hold:*

*a) $p(n, m) - 1 \geq p(n+1, m)$*

*b) $p(n, m) \geq p(n, k) + p(n, m-k)$ for $0 \leq k \leq m$.*

*c) $p(0, m) - 1 \geq 2p(0, m/2)$ for $m$ is even.*

*d) $p(n, m) - 1 \geq p(0, m - 2^{n+1} + 1)$*

*Proof.* a) Since $-1 \leq s(n,m) - s(n+1,m) \leq 0$ because Lemma 6, part a), we have
$$p(n,m) - 1 - p(n+1,m) = 2^{n+3} - 2^{n+2} - 1 - s(n,m) + s(n+1,m)$$
$$= 2^{n+2} - 1 + s(n,m) - s(n+1,m)$$
$$\geq 2^{n+2} - 2$$
$$\geq 0$$

b) Using sub-additivity of $s(n,m)$ (Lemma 5) and the fact that $2^{n+2} - 4 \geq 0$ for $n\mathbb{N}$, we have
$$p(n,m) - p(n,k) - p(n,m-k) \geq s(n,m-k) + s(n,k) - s(n,m) + 2^{n+2} - 4$$
$$\geq s(n,m-k) + s(n,k) - s(n,m)$$
$$\geq 0$$

c) We have
$$p(0,m) - 1 - 2p\left(0, \frac{m}{2}\right) \geq 2s\left(0, \frac{m}{2}\right) - 1 - s(0,m)$$
$$\geq 0,$$
where the last inequality above follows from Lemma 6, part c).

d) We have
$$p(n,m) - 1 - p(0, m - 2^{n+1} + 1) \geq s(0, m - 2^{n+1} + 1) + 1 - s(n,m)$$
$$\geq 0,$$
where last inequality above follows from Lemmas 5 and 6, part b), namely
$$s(n,m) \leq s(n, m - 2^{n+1} + 1 + 2^{n+1} - 1) \leq s(n, m - 2^{n+1} + 1) + s(n, 2^{n+1} - 1)$$
$$\leq s(0, m - 2^{n+1} + 1) + 1$$

This completes the proof. $\square$

We now have presented all lemmas needed to prove the following theorem.

**Theorem 9.** *Let $m, n \in \mathbb{N}$ and assume $m \geq 2^{n+1} - 1$. Then $2^{p(n,m)}\beta_{n,m}$ is an integer, i.e.*
$$-\nu(\beta_{n,m}) \leq p(n,m) \tag{20}$$

*Proof.* From (8) we have
$$2^{p(n,m)}\beta_{n,m} = 2^{p(n,m)-1}\left[\beta_{n+1,m} - \sum_{k=2^{n+1}-1}^{m-2^{n+1}+1} \beta_{n,k}\beta_{n,m-k} - \beta_{0,m-2^{n+1}+1}\right]$$
$$= 2^{p(n,m)-1}\beta_{n+1,m} - \sum_{k=2^{n+1}-1}^{m-2^{n+1}+1} 2^{p(n,m)-1}\beta_{n,k}\beta_{n,m-k} - 2^{p(n,m)-1}\beta_{0,m-2^{n+1}+1} \tag{21}$$

It suffices to show that each term on the right-hand side of (21) is an integer by induction on $m$, which we will do so using properties of $p(n,m)$ established in Lemma 8. Assume that the values of $\beta_{n,m}$ are arranged in a two-dimensional array where the rows are indexed by $n$ and the columns indexed by $m$. Since $\beta_{n,m} = 0$ for $n \geq 1$ and $1 \leq m \leq 2^{n+1} - 2$, we shall call the values in this range trivial and those outside this range, i.e., $m \geq 2^{n+1} - 1$, nontrivial. It follows that each column has at most a finite number of non-trivial entries.

Therefore, we shall apply induction by moving upwards along the non-trivial values in each column from left to right as employed by Ewing and Schober in their induction arguments in [9]. We first establish the base case. Assume $n = 0$ and $m = 1$. Since $\beta_{0,1} = -1/2$ and $p(0,1) = 1$, it is clear that $2^{p(0,1)}\beta_{0,1} = -1$ is an integer.

Next, given any two integers $m \geq 1$ and $n \geq 0$ with $m \geq 2^{n+1} - 1$, we assume inductively that $2^{p(j,k)}\beta_{j,k}$ is an integer for $0 \leq j \leq n$ and $2^{n+1} - 1 \leq k \leq m - 1$; moreover, we also assume that $2^{p(j,m)}\beta_{j,m}$ is an integer for $j \geq n + 1$. Let us consider the first term $2^{p(n,m)-1}\beta_{n+1,m}$ on the right-hand side of (21). Since $p(n, m) - 1 \geq p(n + 1, m)$ (due to part a) in Lemma 8) and $2^{p(n+1,m)}\beta_{n+1,m}$ is an integer by the assumption, it follows that $2^{p(n,m)-1}\beta_{n+1,m}$ is an integer.

Next, we rewrite the summation in (21) according to whether $m$ is even or odd by using Lemma 7. If $m$ is odd, then
$$\sum_{k=2^{n+1}-1}^{m-2^{n+1}+1} 2^{p(n,m)-1}\beta_{n,k}\beta_{n,m-k} = \sum_{k=2^{n+1}-1}^{(m-1)/2} 2^{p(n,m)}\beta_{n,k}\beta_{n,m-k}$$
Since $p(n, m) \geq p(n, k) + p(n, m - k)$ for $0 \leq k \leq m$ from part b) of Lemma 8 and
$$(2^{p(n,k)}\beta_{n,k})(2^{p(n,m-k)}\beta_{n,m-k})$$
is an integer by the assumption, it follows that each term $2^{p(n,m)-1}\beta_{n,k}\beta_{n,m-k}$ in the summation must be an integer. On the other hand, if $m$ is even, then
$$\sum_{k=2^{n+1}-1}^{m-2^{n+1}+1} 2^{p(n,m)-1}\beta_{n,k}\beta_{n,m-k} = \sum_{k=2^{n+1}-1}^{m/2-1} 2^{p(n,m)}\beta_{n,k}\beta_{n,m-k} + 2^{p(n,m)-1}\left(\beta_{n,m/2}\right)^2$$
By the same argument as before, we have that $2^{p(n,m)}\beta_{n,k}\beta_{n,m-k}$ is an integer. Morever, since $p(n, m) - 1 \geq 2p(n, m/2)$ (due to part c) in Lemma 8) and $2^{p(n,m/2)}\beta_{n,m/2}$ is an integer by the assumption, it follows that $2^{p(n,m)-1}\left(\beta_{n,m/2}\right)^2$ must also be an integer. Thus, each term $2^{p(n,m)-1}\beta_{n,k}\beta_{n,m-k}$ in the summation must also be an integer.

As for the last term $2^{p(n,m)-1}\beta_{0,m-2^{n+1}+1}$ in (21), we know from part d) of Lemma 8 that $p(n, m) - 1 \geq p(0, m - 2^{n+1} + 1)$. Since $2^{p(0,m-2^{n+1}-1)}\beta_{0,m-2^{n+1}+1}$ is an integer by the assumption, it follows by the same reasoning that $2^{p(n,m)-1}\beta_{0,m-2^{n+1}+1}$ must be an integer. This finishes the proof of Theorem 9. □

## 3. Special Values of $\beta_{n,m}$

In this section we derive recurrences for special values of $\beta_{n,m}$ where $m$ is restricted to a certain interval. Recall that $\beta_{n,m} = 0$ for $1 \leq m \leq 2^{n+1} - 2$. We therefore begin with an unpublished result by Malik Ahmed and one of the authors regarding $\beta_{n,m}$ in the interval $2^{n+1} - 1 \leq m \leq 2^{n+2} - 3$.

**Theorem 10** (Ahmed-Nguyen). *Let $n \in \mathbb{N}$ and $m$ be a positive integer satisfying $2^{n+1} - 1 \leq m \leq 2^{n+2} - 3$. Then for all $p \in \mathbb{N}$, we have*

$$\beta_{n,m} = \beta_{n+p,m+2^{n+1}(2^p-1)} = -\frac{1}{2}\beta_{0,m-2^{n+1}+1} \tag{22}$$

*Proof.* It follows from (8) that
$$\beta_{n,m} = -\frac{1}{2}\beta_{0,m-2^{n+1}+1}. \tag{23}$$

Next, set
$$n' = n + p, m' = m + 2^{n+1}(2^p - 1).$$

Then
$$m' - 2^{n'+1} + 1 = m - 2^{n+1} + 1$$

which proves
$$\beta_{n,m} = \beta_{n',m'} \tag{24}$$

as desired. □

As a corollary of Theorem 10, we establish a special case of (9).

**Corollary 11.** *Let $n$ be a positive integer and $m$ a positive integer satisfying $2^{n+1} \leq m \leq 2^{n+2} - 3$. Then $2^{2m+2-2^{n+2}}\beta_{n,m}$ is an integer.*

*Proof.* We know from (9) that

$$2^{2(m-2^{n+1}+1)+3-2^2}\beta_{0,m-2^{n+1}+1} = 2^{2m+1-2^{n+2}}\beta_{0,m-2^{n+1}+1}$$

is an integer. It follows from Theorem 10 that

$$2^{2m+2-2^{n+2}}\beta_{n,m} = 2^{2m+2-2^{n+2}}\left(-\frac{1}{2}\beta_{0,m-2^{n+1}+1}\right) = -2^{2m+1-2^{n+2}}\beta_{0,m-2^{n+1}+1} \qquad (25)$$

must also be an integer. □

NOTE: Observe that the corollary above fails for $m = 2^{n+1} - 1$. By Theorem 10 we have $\beta_{n,2^{n+1}-1} = -\frac{1}{2}\beta_{0,0}$. But (9) doesn't apply to $\beta_{0,0} = 1$.

We next focus on deriving recurrences for special values of $\beta_{n,m}$ located at certain positions for $m$ between $2^{n+2} - 2$ and $2^{n+2} + 6$.

**Lemma 12.** *Let $n \in \mathbb{N}$. Then*

$$\beta_{n,2^{n+2}-2} = -\frac{1}{2}\left(\beta_{0,2^{n+1}-1} + \frac{1}{4}\right) \qquad (26)$$

$$\beta_{n,2^{n+2}-1} = -\frac{1}{2}\left(\beta_{0,2^{n+1}} + \frac{1}{4}\right) \qquad (27)$$

*Proof.* Recall that $\beta_{n,m} = -\frac{1}{2}\beta_{0,m-2^{n+1}+1}$ for $n \geq 0$ and $2^{n+1} - 1 \leq m \leq 2^{n+2} - 3$. We have

$$\beta_{n,2^{n+2}-2} = \frac{1}{2}\left[\beta_{n+1,2^{n+2}-2} - \sum_{k=2^{n+1}-1}^{2^{n+1}-1}\beta_{n,k}\beta_{n,2^{n+2}-2-k} - \beta_{0,2^{n+1}-1}\right]$$

$$= \frac{1}{2}\left[0 - \beta_{n,2^{n+1}-1}^2 - \beta_{0,2^{n+1}-1}\right]$$

$$= \frac{1}{2}\left[-\frac{1}{4}\beta_{0,0}^2 - \beta_{0,2^{n+1}-1}\right]$$

$$= -\frac{1}{2}\left[\beta_{0,2^{n+1}-1} + \frac{1}{4}\right]$$

and

$$\beta_{n,2^{n+2}-1} = \frac{1}{2}\left[\beta_{n+1,2^{n+2}-1} - \sum_{k=2^{n+1}-1}^{2^{n+1}}\beta_{n,k}\beta_{n,m-k} - \beta_{0,2^{n+1}}\right] \qquad (28)$$

$$= \frac{1}{2}\left[\beta_{n+1,2^{n+2}-1} - 2\left(\beta_{n,2^{n+1}-1}\beta_{n,2^{n+1}}\right) - \beta_{0,2^{n+1}}\right] \qquad (29)$$

$$= \frac{1}{2}\left[\left(-\frac{1}{2}\beta_{0,0}\right) - 2\left(\left(-\frac{1}{2}\beta_{0,0}\right)\left(-\frac{1}{2}\beta_{0,1}\right)\right) - \beta_{0,2^{n+1}}\right] \qquad (30)$$

$$= -\frac{1}{2}\left[\beta_{0,2^{n+1}} + \frac{1}{4}\right] \qquad (31)$$

□

In the case where $m = 2^{n+2}$, we find that $\beta_{n,m}$ is constant.

**Lemma 13.** *Let $n$ be a positive integer. Then $\beta_{n,2^{n+2}} = 1/16$.*

*Proof.* Recall from Theorem 10 that $\beta_{n,m} = -\frac{1}{2}\beta_{0,m-2^{n+1}+1}$ for $n \geq 0$ and $2^{n+1} - 1 \leq m \leq 2^{n+2} - 3$. Moreover, recall that $b_{2^{n+1}} = 0$ for $n \geq 1$ (Jungreis [12]). Using these results, we have

$$\beta_{n,2^{n+2}} = \frac{1}{2}\left[\beta_{n+1,2^{n+2}} - \sum_{k=2^{n+1}-1}^{2^{n+1}+1} \beta_{n,k}\beta_{n,2^{n+2}-k} - \beta_{0,2^{n+1}+1}\right]$$

$$= \frac{1}{2}\left[-\frac{1}{2}\beta_{0,1} - 2\beta_{n,2^{n+1}-1}\beta_{n,2^{n+1}+1} - \beta_{n,2^{n+1}}^2 - b_{0,2^{n+1}}\right]$$

$$= \frac{1}{2}\left[-\frac{1}{2}\beta_{0,1} - \frac{1}{2}\beta_{0,0}\beta_{0,2} - \frac{1}{4}\beta_{0,1}^2 - 0\right]$$

$$= \frac{1}{2}\left[-\frac{1}{2}(-1/2) - \frac{1}{2}(1)(1/8) - \frac{1}{4}(-1/2)^2\right]$$

$$= 1/16$$

□

We end this section by considering three other special cases.

**Lemma 14.** *Let $n \in \mathbb{N}$. Then*

a) $\beta_{n,2^{n+2}+2} = -\frac{1}{2}\beta_{0,2^{n+1}+3}$ *for $n \geq 2$.*

b) $\beta_{n,2^{n+2}+4} = -\frac{1}{2}\beta_{0,2^{n+1}+5}$ *for $n \geq 2$.*

c) $\beta_{n,2^{n+2}+6} = -\frac{1}{2}\beta_{0,2^{n+1}+7}$ *for $n \geq 3$.*

*Proof.* We have

$$\beta_{n,2^{n+2}+2} = \frac{1}{2}\left[\beta_{n+1,2^{n+2}+2} - \sum_{k=2^{n+1}-1}^{2^{n+1}+3} \beta_{n,k}\beta_{n,2^{n+2}+2-k} - \beta_{0,2^{n+1}+3}\right]$$

$$= \frac{1}{2}\left[\beta_{n+1,2^{n+2}+2} - 2\sum_{k=2^{n+1}-1}^{2^{n+1}} \beta_{n,k}\beta_{n,2^{n+2}+2-k} - \beta_{n,2^{n+1}+1}^2 - \beta_{0,2^{n+1}+3}\right]$$

$$= \frac{1}{2}\left[-\frac{1}{2}\beta_{0,3} - \frac{1}{2}\sum_{j=0}^{1}\beta_{0,j}\beta_{0,4-j} - \frac{1}{4}\beta_{0,2}^2 - \beta_{0,2^{n+1}+3}\right]$$

$$= \frac{1}{2}\left[-\frac{1}{2}(-1/4) - \frac{1}{2}(\beta_{0,0}\beta_{0,4} + \beta_{0,1}\beta_{0,3}) - \frac{1}{4}(1/8)^2 - \beta_{0,2^{n+1}+3}\right]$$

$$= \frac{1}{2}\left[-\frac{1}{2}(-1/4) - \frac{1}{2}[(1)(15/128) + (-1/2)(-1/4)] - \frac{1}{4}(1/8)^2 - \beta_{0,2^{n+1}+3}\right]$$

$$= -\frac{1}{2}\beta_{0,2^{n+1}+3}$$

This proves part a). Parts b) and c) can be proven in a similar manner. □

## 4. New Area Approximations

In this section we describe a parallel processing algorithm to compute the values of $\beta_{n,m}$ and present new upper bounds for the area of $M$ that were calculated using this algorithm. Assume as before that the values of $\beta_{n,m}$ are arranged in a two-dimensional array with the rows indexed by $n$ and columns indexed by $m$. We recall Ewing and Schober's backwards algorithm for computing the non-trivial values of $\beta_{n,m}$ recursively one at a time by moving upwards along each column from left to right as described in our induction proof of Theorem 9. Thus, the order of computation would be:

$$\beta_{0,1}, \beta_{0,2}, \beta_{1,3}, \beta_{0,3}, \beta_{1,4}, \beta_{0,4}, ...$$

Our new method is as follows: we calculate values of $\beta_{n,m}$ across multiple columns simultaneously in a parallel fashion while moving up along them as before until we reach a critical row near the top where from this point on, all remaining column values must be computed one at a time. This is then repeated for the next set of columns, etc.

To illustrate our method, consider for example the calculation of $\beta_{1,7}$ and $\beta_{1,8}$ in row $n=1$ using the backward recursion formula (8):

$$\beta_{1,7} = \frac{1}{2}\left[\beta_{2,7} - \sum_{k=3}^{4}\beta_{1,k}\beta_{1,7-k} - \beta_{0,4}\right] = \frac{1}{2}\left[\beta_{2,7} - 2\beta_{1,3}\beta_{1,4} - \beta_{0,4}\right]$$

$$\beta_{1,8} = \frac{1}{2}\left[\beta_{1,8} - \sum_{k=3}^{4}\beta_{1,k}\beta_{1,8-k} - \beta_{0,5}\right] = \frac{1}{2}\left[\beta_{2,8} - 2\beta_{1,3}\beta_{1,5} - \beta_{1,4}^2 - \beta_{0,4}\right]$$

Observe that these two values do not depend on each other and can be computed independently in parallel. However, this is not the case for $\beta_{0,7}$ and $\beta_{0,8}$ in the top row ($n=0$) where the latter depends on the former:

$$\beta_{0,7} = \frac{1}{2}\left[\beta_{1,7} - \sum_{k=1}^{6}\beta_{0,k}\beta_{0,7-k} - \beta_{0,6}\right] = \frac{1}{2}\left[\beta_{1,7} - 2\beta_{0,1}\beta_{0,6} - 2\beta_{0,2}\beta_{0,5} - 2\beta_{0,3}\beta_{0,4} - \beta_{0,4}\right]$$

$$\beta_{0,8} = \frac{1}{2}\left[\beta_{1,8} - \sum_{k=1}^{7}\beta_{0,k}\beta_{0,8-k} - \beta_{0,7}\right] = \frac{1}{2}\left[\beta_{1,8} - 2\beta_{0,1}\beta_{0,7} - 2\beta_{0,2}\beta_{0,6} - 2\beta_{0,3}\beta_{0,5} - \beta_{0,4}^2 - \beta_{0,7}\right]$$

In general, the values $\beta_{n,m}$, $\beta_{n,m+1}$ and $\beta_{n,m+2}$ in three consecutive columns can be calculated in parallel as long as $n \geq 1$. This is because $\beta_{n,m+1}$ depends only on the values $\beta_{n,k}$ in row $n$, where $k = 3, 4, ..., m-2$, which are prior to $\beta_{n,m}$. Similarly, $\beta_{n,m+2}$ depends only on $\beta_{n,k}$ where $k = 3, 4, ..., m-1$. Since the number of non-zero values in each column increases as $m$ increases, this parallel algorithm becomes more effective and asymptotically three times as fast in comparison to that of calculating $\beta_{n,m}$ one at a time. Moreover, this approach can be extended to calculate the values $\beta_{n,m}, \beta_{n,m+1}, ..., \beta_{n,m+6}$ in seven consecutive columns simultaneously as long as $n \geq 2$. More generally, if $n \geq N$, then up to $2^{N+1} - 1$ columns can be computed in parallel.

We were able to use this parallel algorithm to calculate the first five million terms of $b_m$ and obtain a new upper bound of $A_{5,000,000} \approx 1.68288$ for the area of the Mandelbrot set. This algorithm was implemented using the programming language C++ and message passing interface Open MPI. In particular, we calculated the values of $\beta_{n,m}$ across four columns in parallel for $n \geq 2$ beginning with the first group of columns $\beta_{n,8}$, $\beta_{n,9}$, $\beta_{n,10}$, $\beta_{n,11}$ (we initialized columns $\beta_{n,0}, ..., \beta_{n,7}$ with their known values). Our code was executed on a Linux cluster and required four processors (1.05 Ghz AMD Opteron 2352 quad-core processors) to execute it since each column was computed using a different processor. After computing its column of values, each processor would pass these values to the other three processors before calculating to its next designated column. Thus, each processor was required to store all values of $\beta_{m,n}$ (generated from all four processors) separately in its own RAM in order to compute its next column. This parallel approach improved the performance of our implementation significantly; asymptotically, run-time was decreased by a factor of four in comparison to using a single processor, but at the cost of quadrupling our RAM memory requirements. The only computational cost to our algorithm involves having each processor pass its values to the other three processors. Since the number of non-zero values for $\beta_{n,m}$ in each column grows on the order of $\log_2 m$, the computational cost in passing these values is insignificant in comparison to the cost of computing $\beta_{n,m}$ itself using (8), whose summation term grows on the order of $m$ since $2^{n+1} - 1 \leq m \leq 2^{n+2} - 2$.

Table 1 gives values for the approximations $A_N$, where $N$ ranges from 500,000 to 5 million in increments of 500,000, based on our computed values of $\beta_{n,m}$, and thus $b_m = \beta_{0,m+1}$. These values were computed in batches over a five-month period between August-December, 2014, although the actual total run-time was approximately 3 months. Table 2 gives the reader a sense of the run-time required to compute $b_m$ in batches of 500,000 starting at $m = 2,500,000$.

To estimate the error in our upper bounds, we use Ewing and Schober's [9] analysis of their calculation of $\beta_{n,m}$ using (8) and double-precision floating-point arithmetic. First, they considered propagation error due to errors in computing previous coefficients. They argued probabilistically that the propagation error is on the same order of magnitude as machine error so that the computations for $\beta_{n,m}$ are stable. In particular, let

$$\tilde{\beta}_{n,m} = \beta_{n,m} + \epsilon_{n,m}, \tag{32}$$

where $\beta_{n,m}$ is the true value, $\tilde{\beta}_{n,m}$ is the calculated value, and $\epsilon_{n,m}$ is the corresponding error. Then substituting (32) into (8) shows that the propagation error in computing $\beta_{n,m}$ is on the order of

$$\epsilon_{n,m} = \beta_{n,m} - \tilde{\beta}_{n,m}$$
$$= \frac{1}{2}\left(\beta_{n+1,m} - \tilde{\beta}_{n+1,m} - \sum_{k=2^{n+1}-1}^{m-2^{n+1}+1}(\beta_{n,k}\beta_{n,m-k} - \tilde{\beta}_{n,k}\tilde{\beta}_{n,m-k}) - \beta_{0,m-2^{n+1}+1} + \tilde{\beta}_{0,m-2^{n+1}+1}\right)$$
$$\approx \frac{1}{2}(-\epsilon_{n+1,m} + \epsilon_{0,m-2^{n+1}+1}) + \sum_{k=2^{n+1}-1}^{m-2^{n+1}+1}\beta_{n,k}\epsilon_{n,m-k}, \tag{33}$$

where the quadratic error terms are ignored. Next, assume that $\epsilon_{n,m}$ is uniformly distributed with a small probability of exceeding some threshold value $\epsilon$. Moreover, assume that the sum

$$\sum_{k=2^{n+1}-1}^{m-2^{n+1}+1}|\beta_{n,k}|$$

is bounded, which we verified in computing $A_{5,000,000}$. It follows from the law of large numbers that the error term in (33),

$$\sum_{k=2^{n+1}-1}^{m-2^{n+1}+1}\beta_{n,k}\epsilon_{n,m-k}, \tag{34}$$

which we view as a weighted sum of independent and identically distributed random variables, approaches zero as $m \to \infty$. Thus, (34) is negligible in contributing towards the propagation error in (33). Hence, if all previous errors are bounded by $\epsilon$, then so will the propagation error.

To check the accuracy of our calculations, we compared our calculated values of $b_m$ (in double-precision floating point format) with exact values that are given by closed formulas at certain positions. For example, Levin [14] and Ewing and Schober [9] proved independently that $b_m = 0$ for all $m = (2k+1)2^\nu$, where $k, \nu \in \mathbb{N}$ satisfy $k + 3 \leq 2^\nu$. We confirmed this for our calculated values. Moreover, Ewing and Schober [9] proved that

$$b_m = \begin{cases} \dfrac{-1}{2^{\nu+3}(2^\nu - 1)}\binom{2^\nu - 5/2}{2^\nu - 2}, & m = (2^\nu - 1)2^\nu, \nu \geq 1; \\ \dfrac{3(2^\nu - 6)}{2^{\nu+5}(2^\nu + 1)(2^\nu - 5)}\binom{2^\nu - 3/2}{2^\nu - 1}, & m = (2^\nu + 1)2^\nu, \nu \geq 2; \\ \dfrac{-(214 \cdot 2^{3\nu} - 767 \cdot 2^{2\nu} + 146 \cdot 2^\nu + 452)}{2^{\nu+8}(2^{\nu+1} - 7)(2^{2\nu} - 1)(2^\nu + 2)}\binom{2^\nu - 5/2}{2^\nu - 2}, & m = (2^\nu + 3)2^\nu, \nu \geq 2. \end{cases} \tag{35}$$

A comparison of our calculated values for $b_m$ with the exact values at the positions specified in (35) yielded a maximum error of $5.00034 \cdot 10^{-16}$. Thus, the error in calculating $b_m$ for $m \leq 5,000,000$ is at most $6 \cdot 10^{-16}$ and the error in our upper bound $A_{5,000,000}$ is at most $3 \cdot 10^{-9}$. We note that our calculation of $A_{1,000,000} \approx 1.70393$ agrees with that reported by Chen, Kawahira, Li, and Yuan [5]. In Table 3, we give values for $b_m$ at certain positions where no closed formula is known so that the reader may verify our calculations.

Figure 1 shows a plot of Table 1 that clearly reveals the slow convergence of $A_N$. If the exact value of $A$ lies closer to 1.50659 as computed by pixel counting, then certainly using $A_N$ to closely approximate $A$ is impractical due to the extremely large number of terms required. On the other hand, if the exact value lies closer to 1.68, then this would indicate that the boundary of the Mandelbrot set may have positive area.

TABLE 1. New Upper Bounds for the Area of the Mandelbrot Set

| $N$ (millions) | $A_N$ |
|---|---|
| 0.5 | 1.72 (Ewing-Schober) |
| 1 | 1.70393 |
| 1.5 | 1.69702 |
| 2 | 1.69388 |
| 2.5 | 1.69096 |
| 3 | 1.68895 |
| 3.5 | 1.6874 |
| 4 | 1.68633 |
| 4.5 | 1.68447 |
| 5 | 1.68288 |

TABLE 2. Run-Times for Calculating $b_m$ in Batches of 500,000

| Range of $m$ (millions) | Run-time to compute $b_m$ (days) |
|---|---|
| 2.5-3 | 9 |
| 3-3.5 | 10.8 |
| 3.5-4 | 12.5 |
| 4-4.5 | 14.4 |
| 4.5-5 | 16.2 |

TABLE 3. Some calculated values of $b_m$ at positions where no closed formula is known.

| $m$ | $b_m$ |
|---|---|
| 500,000 | $5.5221313 \cdot 10^{-8}$ |
| 1,000,000 | $-4.713883 \cdot 10^{-8}$ |
| 1,500,000 | $8.4477641 \cdot 10^{-8}$ |
| 2,000,000 | $-6.437866 \cdot 10^{-9}$ |
| 2,500,000 | $1.6594295 \cdot 10^{-8}$ |
| 3,000,000 | $8.150385 \cdot 10^{-9}$ |
| 3,500,000 | $-3.911993 \cdot 10^{-9}$ |
| 4,000,000 | $2.315128 \cdot 10^{-9}$ |
| 4,500,000 | $-8.87746 \cdot 10^{-9}$ |
| 5,000,000 | $8.0532 \cdot 10^{-11}$ |

## 5. Conclusions

In this paper we presented new results which improve on known upper bounds for the area of the Mandelbrot set and 2-adic valuations of the series coefficients $\beta_{n,m}$ given by Ewing and Schober in [9]. Of course, our calculations of the first five million terms of $b_m$ were performed using more powerful computers that those available to Ewing and Schober two decades ago. Therefore, it would be interesting to find out in the next two decades what improvements can be made to our results by using computers that will be even more powerful, unless we are fortunate enough to see the exact area calculated before then.

## 6. Acknowledgments

The authors wish to thank the reviewer for providing us with feedback that greatly improved this paper and bringing to our attention important references that allowed us to give a more accurate survey of previous results. We also wish to thank Mark Sedlock and Jason Faulcon for their technical support of the Computing Cluster at Rowan University on which our computations were performed.

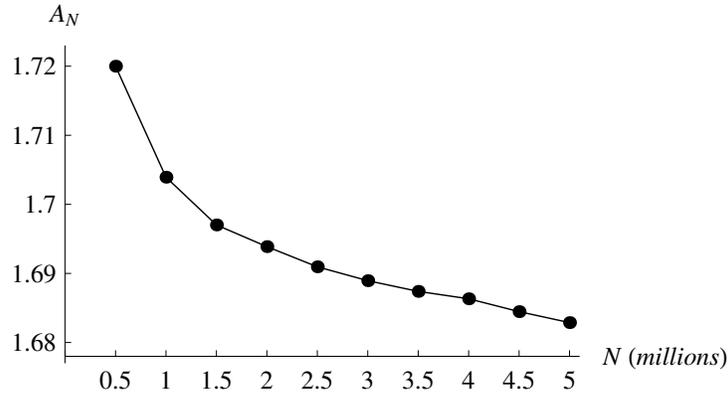

Figure 1. Plot of $A_N$

Department of Mathematics, Rowan University, Glassboro, NJ 08028.
*E-mail address*: `bittne12@students.rowan.edu`

*E-mail address*: `cheong94@students.rowan.edu`

*E-mail address*: `gatesd78@students.rowan.edu`

*E-mail address*: `nguyen@rowan.edu`